\documentclass[a4paper,12pt]{article}
\usepackage{amssymb,amsmath,amsthm,latexsym}
\usepackage[croatian,english]{babel}

\usepackage{color}

\newcommand{\normal}{\color{black}}

\usepackage[a4paper,margin=1in]{geometry} 
\usepackage{verbatim}
\usepackage{epsf,graphicx}

\theoremstyle{plain}
\newtheorem{thm}{Theorem}[section]

\newtheorem{corollary}[thm]{Corollary}
\newtheorem{lemma}[thm]{Lemma}

\theoremstyle{definition}

\numberwithin{equation}{section}

\newcommand{\R}{{\mathbb R}}
\newcommand{\N}{{\mathbb N}}

\newcommand{\E}{\mathbb{E}\,}

\def\P{{\mathbb P}}

\begin{document}

\title{\Large\bfseries Distribution of suprema for generalized risk processes}

\author{\itshape
    Ivana Ge\v{c}ek Tu{\dj}en\\
    \\
    Department of Mathematics\\
    University of Zagreb, Zagreb, Croatia\\
    Email: igecek@math.hr
   \normal
}

\date{}

\maketitle

\begin{abstract}
\noindent
We study a generalized risk process $X(t)=Y(t)-C(t)$, $t\in[0,\tau]$, where $Y$ is a L\'evy process, $C$ an independent subordinator and $\tau$ an independent exponential time. Dropping the standard assumptions on the finite expectations of the processes $Y$ and $C$ and the net profit condition, we derive a Pollaczek-Khinchine type formula for the supremum of the dual process $\widehat{X}=-X$ on $[0,\tau]$ which generalizes the results obtained in \cite{HPSV1}. We also discuss which assumptions are necessary for deriving this formula, specially from the point of view of the ladder process.


\end{abstract}

\section{Introduction}\label{sec-1}
The basic risk model, known as the \emph{Cram\'{e}r-Lundberg model}, has been revisited many times in the risk theory. It is a model that is based on the risk process $(R(t)~:~t\geq 0)$ such that
$$R(t)=ct-\sum_{i=1}^{N(t)} Y_i~,~~t\geq 0~,$$ where $c>0$ represents the premium rate (we assume that there are the incoming premiums which arrive from the policy holders),
$(Y_i~:~i\in \N)$ is an i.i.d. sequence of nonnegative random variables with common distribution $F$ (which usually represent the policy holders' claims) and $(N(t)~:~t\geq 0)$ a homogeneous Poisson process of rate $\lambda>0$, independent of $(Y_i~:~i\in \N)$.
One of the main questions that is observed in this model is the question of the \emph{ruin probability}, given some initial capital
$u>0$, i.e.
$$\vartheta(u)=\P(u+R(t)<0~,~~\textrm{for some $t>0$})~.$$
Sometimes it will be easier to work with the \emph{survival probability}, so we also define
$$\theta(u)=1-\vartheta(u)~.$$
This model was generalized by a few authors and the generalization usually goes two ways. First we can allow additional uncertainties in income (premiums), which means that we can add a perturbation in the basic model with the compound Poisson process. The perturbation is usually modelled by a diffusion and methods from martingale theory or renewal theory are used to gain some results. This approach was used for example by \emph{Dufresne} and \emph{Gerber} in \cite{DG}, for the standard risk process which is perturbed by the standard Brownian motion. But if we model the perturbation by, for example, gamma process then fluctuation theory for L\'evy processes is very natural and applicable to use, especially since many risk processes are special cases of spectrally negative L\'evy processes. This was first done by \emph{Furrer} (in \cite{Furr}) who used \emph{Zolotarev's} result (for details we refer to \cite{Zol}) which establishes a connection between the distribution of the infimum of the $\alpha$-stable L\'evy process and its characteristic exponent. We will also use results from the fluctuation theory for spectrally negative L\'evy processes to gain our results in this paper.\\
\\
Another possible generalization can be made in the claim process (in generalization mentioned before it was still modelled by the compound Poisson process). This was done by \emph{Dufresne}, \emph{Gerber} and \emph{Shiu} in \cite{DGS} where the claim process was modelled by a gamma process. This model was further generalized by Yang and Zhang in \cite{YZ} where the gamma process was perturbed by the Brownian motion and the authors used \emph{Furrer's} approach in decomposing the supremum of the dual of the risk process at the modified ladder epochs (which were first introduced by \emph{Schmidli}, see \cite{Schm} for details). This model was then generalizd by \emph{Huzak}, \emph{Perman}, \emph{\v Siki\' c} and \emph{Vondra\v cek} in \cite{HPSV1} where the authors considered the generalized risk process $$X(t)=ct-C(t)+Z(t)~,~~t\geq 0~,$$ for
$(C(t)~:~t\geq 0)$ a subordinator and $(Z(t)~:~t\geq 0)$ an independent spectrally negative L\'{e}vy
process. They obtained a \emph{Pollaczek-Khinchine} type result, i.e.
$$\theta(u)=(1-\rho)\cdot\sum_{n=0}^{\infty}\rho^n(G^{(n+1)\ast}\ast
H^{n\ast})(u)~,~~u\geq 0~,$$ where $G$ is the distribution function of the supremum of the dual process $\widehat{X}=-X$ before time $\sigma$, the first time the new supremum is reached because of the jump of the process $C$, and $H$ the integrated tail of the distribution of the jumps only related to the process $C$.\\
\\
In this paper we will also consider a generalized risk process, in the sense that we will take a generalized claim process and also allow a perturbation in our model. More precisely, we will take a L\'{e}vy process $Y$ and an independent subordinator $C$ and consider the L\'{e}vy process
\begin{equation}\label{1:eq.1}
X(t):=Y(t)-C(t)~,~~t\in[0,\tau]~,
\end{equation}
where $\tau$ is some independent exponential time, $\tau=\tau(q)\sim Exp(q)$, $q>0$.\\
\\
In the basic risk model one usually assumes that the expectations of the processes involved are finite and that
$$\E R(1)=c-\lambda\E Y(1)>0~,$$
which is called the \emph{net profit condition}. Given that condition,
$\lim_{t\to\infty} R(t)=+\infty$, i.e. the ruin probability is not equal to $1$. In case that $c \leq \lambda \mu$, we have $\E R(1)\leq 0$. If this expectation is zero, then $-\infty = \liminf R(t) < \limsup R(t) = +\infty$ and if $\E R(1)<0$, then $\lim_{t\to\infty}
R(t)=-\infty$, which means that in both cases $\vartheta(u)=1$.\\
\\
In this paper we will drop the standard assumptions on finite expectation of the underlying processes and net profit condition and we will show that in this generalized setting we can again derive a \emph{Pollaczek-Khinchine} type formula. More precisely, for the dual process $\widehat{X}=-X$ and $\widehat{S}(t)=\sup_{0\leq s\leq t}\widehat{X}(s)$ we will define
\begin{equation}\label{1:eq.2}
\sigma=\inf\{t>0~:~\Delta
C(t)>\widehat{S}(t-)-\widehat{X}(t-)\}
\end{equation}
as the first time when the new supremum of the dual process $\widehat{X}$ is reached because of the jump of the subordinator $C$. Following the approach used in \cite{Furr}, \cite{Schm} and \cite{HPSV1}, we will decompose our dual process at the successive times $\sigma_i$, $i\in\N$, and
we will derive a \emph{Pollaczek-Khinchine} type formula for the supremum of the dual process $\widehat{X}$ on $[0,\tau]$ in this generalized setting. In order to do this, we will need some basic assumptions without which decomposition would not make sense, i.e. we will assume that $\P(\sigma>0)=1$ and $\lim_{t\to\infty} X(t)=\infty$, a.s. Under these assumptions we will show that
\begin{equation}\label{1:eq.3}\P(\widehat{S}(\tau)\leq
x)=(1-\rho_\tau)\sum_{n=0}^\infty\rho_\tau^n(G_\tau^{(n+1)\ast}\ast
H_\tau^{n\ast})(x)~,~~x\geq 0~,\end{equation} where $G_\tau$ is the distribution function of $\widehat{S}((\sigma\wedge\tau)-)$,
\begin{equation*}\rho_\tau=\P(\sigma\leq\tau)=\frac{\kappa(q,0)}{q}\P(\sigma>\tau)\int_0^{\infty}\nu(u,\infty)\Upsilon^q(du)\end{equation*}
and
\begin{equation*}H_\tau(x)=\frac{\int_0^x \nu(u,\infty)\Upsilon^q(du)}{\int_0^\infty
\nu(u,\infty)\Upsilon^q(du)}~,\end{equation*}
for $\nu$ being the L\'{e}vy measure of the subordinator $C$,
\begin{equation*}\Upsilon^q(x)=\int_0^\infty \exp\{-qL^{-1}(t)\}\P(H(t)\leq
x)\,dt~,~~x\geq 0~,\end{equation*} the modified renewal function (where $(L^{-1},H)$ is the ladder process, precise definition will be given in Section 2.) and $\kappa$ comes from its Laplace transformation, \begin{equation*}\lambda\int_0^\infty
e^{-qx}\Upsilon^q(x)\,dx=\frac{1}{\kappa(q,\lambda)}~,~~\lambda>0~.\end{equation*}
As a consequence, if we consider a standard generalized risk process
\begin{equation}\label{1:eq.4}
X(t)=ct-C(t)+Z(t)~,~~t\in[0,\tau]~,\end{equation} where
$C$ is a subordinator and $Z$ a spectrally negative L\'{e}vy
process (and we are assuming that these processes are independent) and if we take the standard assumptions on finite expectation (i.e. $\E[C(1)]<\infty$ and, without loss of generality, $\E[Z(1)]=0$) and net profit condition (i.e. $\E[C(1)]<c$), the next more specific result can be obtained:
\begin{equation}\label{1:eq5}\P(\widehat{S}(\tau)\leq x)=(1-p_\tau)\sum_{n=0}^{\infty}
{p_\tau}^n({\widetilde{G}_\tau}^{(n+1)*}*{\widetilde{H}_\tau}^{n*})(x)~,~x\geq 0~,\end{equation}
where $\widetilde{G}_\tau$ is the distribution function of $\widehat{S}((\sigma\wedge\tau)-)$, $p_\tau=\P(\sigma\leq\tau)$ and $$\widetilde{H}_\tau(x)=\frac{1}{\int_0^\infty e^{-\phi(q) u}\nu(u,\infty)
du}\int_0^x e^{-\phi(q) u}\nu(u,\infty)\,du~,$$ the integrated tail of the distribution of the jumps of the subordinator $C$. Here $\nu$ is the L\'evy measure of the subordinator $C$ and $\phi=\phi_X=\psi_X^{-1}$ is the inverse of the Laplace exponent of the process $X$.\\
\\
Finally, in the Section 3 we will discuss our results from the point of view of the ladder height process, explaining the role of the net profit condition and other assumptions for the results on the supremum of (the dual of) the generalized risk process.

\section{Distribution of the supremum for generalized risk process stopped at independent exponential time}\label{sec-2}
Let $Y=(Y(t)~:~t\geq 0)$ be a L\'evy process with corresponding L\'evy measure $\Pi$ (i.e. $\Pi$ is a measure on $\R\setminus\{0\}$ such that $\int_{\R}(x^2\wedge 1)\Pi(dx)<\infty)$) and $C=(C(t)~:~t\geq 0)$ a subordinator independent of $Y$ with corresponding L\'evy measure $\nu$ (i.e. $\nu$ is a measure on $(0,\infty)$ such that $\int_{(0,\infty)}(x\wedge 1)\nu(dx)<\infty$). We define a process $X$ on $[0,\tau]$ by
\begin{equation}\label{2:eq.1}
X(t)=Y(t)-C(t)~,~~t\in[0,\tau]~,
\end{equation}
where $\tau=\tau(q)$, $q>0$, is an independent exponential time with parameter $q$. Let us denote by $\mathcal{F}=(\mathcal{F}(t)~:~t\geq 0)$ the standard augmentation of $\mathcal{F}^0=(\mathcal{F}^0(t)~:~t\geq 0)$, where $\mathcal{F}^0(t)=\sigma(Y_s,C_s~:~0\leq s\leq t)$, and let us notice that $X$ is a L\'evy process with respect to the filtration $\mathcal{F}$.\\
\\
Let us introduce notation for supremums which we will soon need:
\begin{equation*}
\widehat{S}(t)=\sup_{0\leq s \leq
t} \widehat{X}(s)~~\textrm{and}~~\widehat{S}(\infty)=\sup_{0\leq s < \infty}
\widehat{X}(s)~.
\end{equation*}
For the original process $X$ we can define infimums,
$I(t):=\inf_{0\leq s\leq t} X(s)$ and $I(\infty):=\inf_{0\leq
s<\infty} X(s)$.\\
\\
Now we observe the dual process $\widehat{X}=-X$ and we would like to decompose it at some specific stopping times, which we will call \emph{modified ladder times}. These are the times in which process $\widehat{X}$ reaches the new supremum because of the jump of the subordinator $C$, i.e. times $t$ when $\Delta C(t)>\widehat{S}(t-)-\widehat{X}(t-)$. So we define the first of such times as
\begin{equation}\label{2:eq.2}
\sigma=\inf\{t>0~:~\Delta
C(t)>\widehat{S}(t-)-\widehat{X}(t-)\}~.
\end{equation}
Without the assumptions on the finite expectations of the underlying processes and net profit condition, a natural question arises: can we decompose our supremum at these times? Namely, the set of times when the new supremum of the dual process is reached because of the jump of the subordinator $C$ needs not to be discrete. More precisely, using the Blumenthal's law, we know that $\P(\sigma=0)=0$ or $1$. If $\sigma>0$ a.s., then $\sigma$ is really the first time when the new supremum of the process $\widehat{X}$ is reached because of the jump of the subordinator $C$ and if $\sigma=0$ a.s., then number of these times is infinite , i.e. $0$ is an accumulation point of such times. We refer to \cite{SV} for necessary and sufficient conditions for $\sigma>0$ a.s.\\
\\
So, in order to make our decomposition possible, let us assume that
\begin{itemize}\item[(i)] \begin{equation} \label{2:eq.3}\P(\sigma>0)=1~,\end{equation}
\item[(ii)]\begin{equation} \label{2:eq.4} \lim_{t\to\infty} X(t)=\infty~~\textrm{a.s.}\end{equation}
\end{itemize}
Now we are interested in the reflected process $\widehat{S}-\widehat{X}$ and we want to determine the expected time that this process spends in $(0,x)$, $x>0$, until the time $\sigma\wedge\widehat{\tau_y}\wedge\tau$, i.e.  $\E\int_{0}^{\sigma\wedge\widehat{\tau_y}\wedge\tau}
1_{\{\widehat{S}(t)-\widehat{X}(t)\leq x\}} dt$, where
\begin{equation*}
\widehat{\tau_y}=\inf\{t>0~:~\widehat{X}(t)>y\}~,~y>0~,
\end{equation*}
is the first time when the process $\widehat{X}$ enters $(y,\infty)$. This will lead us to the auxiliary result (Lemma 2.1.) which will allow
us to determine the distributions of the crucial variables involved in our problem. We obviously have
\begin{equation*}  \widehat{S}(t-)\leq y ~~\textrm{if and only if}
~~ t\leq \widehat{\tau_y}~.\end{equation*}
 For $x>0$ and $y>0$, using Fubini's theorem, linearity of expectation, continuity in probability of the process $\widehat{X}$ and $S(t)=^d
\widehat{S}(t)-\widehat{X}(t)$, we have
\begin{align}
&\E\int_0^{\tau} 1_{\{\widehat{S}(t)-\widehat{X}(t)\leq x\}}\,dt\nonumber\\
&=\E\int_0^{\infty} qe^{-qs}\left(\int _0^s
1_{\{\widehat{S}(t)-\widehat{X}(t)\leq x\}}\,dt\right)\,ds=\int_0^{\infty}
qe^{-qs}\left(\int_0^s \P(S(t)\leq x)\,dt\right)\,ds\nonumber\\
&=\int_0^{\infty}\P(S(t)\leq x)\left(\int_t^{\infty} qe^{-qs} ds\right)\,dt=\int_0^{\infty} e^{-tq}\P(S(t)\leq x)\,dt\nonumber\\
&=\frac{1}{q}\int_0^{\infty} qe^{-tq}\P(S(t)\leq x)\,dt~.\label{2:eq.5}\end{align}
Let $(L^{-1},H)$ be the \emph{ladder process} associated to the process $X$. It is a L\'evy process and we denote it's \emph{bivariate Laplace exponent} by $\kappa$, i.e.
\begin{equation*}
\exp\{-\kappa(\alpha,\beta)\}=\E[\exp\{-(\alpha L^{-1}(1)+\beta
H(1))\}]~,~~\alpha,\beta>0~.\end{equation*}
We also define \emph{the renewal function} $\Upsilon$ associated to the ladder height process $H$ as
\begin{equation*}
\Upsilon(x)=\int_0^\infty\P(H(t)\leq x)\,dt=\E\int_0^\infty
1_{\{H(t)\leq x\}}\,dt=\E\int_0^\infty 1_{\{S(t)\leq x\}}\,dL(t)~,~x\geq 0~.
\end{equation*}
For $q>0$ we define the function
\begin{equation}\label{2:eq.6}\Upsilon^q(x)=\int_0^\infty \exp\{-qL^{-1}(t)\}\P(H(t)\leq
x)\,dt~,~~x\geq 0~.\end{equation}
The Laplace transform of $\Upsilon^q$ is given by (for details see \cite[p. 172-174]{Ber}):
\begin{equation}\label{2:eq.6a}\lambda\int_0^\infty
e^{-qx}\Upsilon^q(x)\,dx=\frac{1}{\kappa(q,\lambda)}~,~~\lambda>0\end{equation}
and we have
\begin{equation}\label{2:eq.6b}q\int_0^\infty
e^{-qt}\P(\tau_x>t)\,dt=\kappa(q,0)\Upsilon^{q}(x)~.\end{equation}
Now from (\ref{2:eq.5}) it follows that
\begin{equation}\label{2:eq.7}
\E\int_0^\tau 1_{\{\widehat{S}(t)-\widehat{X}(t)\leq x\}}\,
dt=\frac{1}{q}\kappa(q,0)\Upsilon^q(x)~.\end{equation}
Now we look at the expected time that process $\widehat{S}-\widehat{X}$ spends under the level $x$ until the time $\tau$, but after the time
$\sigma\vee\widehat{\tau_y}=\max\{\sigma,\widehat{\tau_y}\}$. We have
\begin{align}
&\E\int_0^{\tau} 1_{\{\widehat{S}(t)-\widehat{X}(t)\leq x\}}
1_{\{t>\sigma\}}1_{\{\widehat{S}(t)>y\}}
\,dt\nonumber\\
&=\E\big[\int_{\sigma\vee
\widehat{\tau_y}}^{\tau} 1_{\{\widehat{S}(t)-\widehat{X}(t)\leq
x\}}
\,dt|\sigma\vee \widehat{\tau_y}\leq\tau\big]\P(\sigma\vee \widehat{\tau_y}\leq\tau)\nonumber\\
&=\P(\sigma\vee \widehat{\tau_y}\leq\tau)\E\int_0^{\tau}1_{\{\widehat{S}(t)-\widehat{X}(t)\leq
x\}}\,dt\nonumber\\
&=\P(\sigma\leq
\tau,\widehat{\tau_y}\leq\tau)\frac{1}{q}\kappa(q,0)\Upsilon^q(x)~.\label{2:eq.8}
\end{align}
Analogusly we get
\begin{equation} \label{2:eq.9}\E\int_0^{\tau} 1_{\{\widehat{S}(t)-\widehat{X}(t)\leq x\}}
1_{\{\widehat{S}(t)>y\}}\,dt=\P(\widehat{\tau_y}\leq\tau)\frac{1}{q}\kappa(q,0)\Upsilon^q(x)~.\end{equation}
Now we subtract (\ref{2:eq.9}) from (\ref{2:eq.8}) and we get
\begin{align}&\E\int_0^{\tau} 1_{\{\widehat{S}(t)-\widehat{X}(t)\leq x\}}
1_{\{t\leq\sigma\}}1_{\{\widehat{S}(t)>y\}}\,dt\nonumber\\
&=\frac{1}{q}\kappa(q,0)\Upsilon^q(x)(\P(\widehat{\tau_y}\leq\tau)-\P(\sigma\leq\tau,\widehat{\tau_y}\leq\tau))\nonumber\\
&=\frac{1}{q}\kappa(q,0)\Upsilon^q(x)\P(\widehat{\tau_y}\leq\tau,\sigma>
\tau) \label{2:eq.10}\end{align}
Similarly we get
\begin{equation}\label{2:eq.11}\E\int_0^{\tau} 1_{\{\widehat{S}(t)-\widehat{X}(t)\leq x\}}
1_{\{t\leq\sigma\}}\,dt=\P(\sigma> \tau)\frac{1}{q}\kappa(q,0)\Upsilon^q(x)~.\end{equation}
Now we subtract (\ref{2:eq.11}) from (\ref{2:eq.10})
\begin{align*}&\E\int_0^{\tau} 1_{\{\widehat{S}(t)-\widehat{X}(t)\leq x\}}
1_{\{t\leq\sigma\}}1_{\{\widehat{S}(t)\leq y\}}\,dt\\
&=\frac{1}{q}\kappa(q,0)\Upsilon^q(x)(\P(\sigma> \tau)-\P(\widehat{\tau_y}\leq\tau,\sigma> \tau))\\
&=\frac{1}{q}\kappa(q,0)\Upsilon^q(x)\P(\sigma>\tau,\widehat{\tau_y}> \tau)~,
\end{align*}
which leads us to the following result.
\begin{lemma} For $x>0$ and $y>0$ we have
\begin{equation}\label{2:eq.12}\E\int_0^{(\sigma\wedge\widehat{\tau_y})\wedge\tau}
1_{\{\widehat{S}(t)-\widehat{X}(t)\leq
x\}}\,dt=\P(\sigma>\tau,\widehat{\tau_y}>
\tau)\frac{1}{q}\kappa(q,0)\Upsilon^q(x)~.\end{equation}
\end{lemma}
\vspace{0.3cm}
We can rewrite (\ref{2:eq.12}) as
\begin{equation*}\E\int_0^{\sigma\wedge\widehat{\tau_y}\wedge\tau}
1_{\{\widehat{S}(t)-\widehat{X}(t)\leq x\}}\,
dt=\frac{\kappa(q,0)}{q}\P(\sigma>\tau,\widehat{\tau_y}>\tau)\int_0^\infty
1_{\{u\leq x\}} \Upsilon^q(dx)~,\end{equation*}
which means that for every nonegative Borel function $f$ we have
\begin{equation}\label{2:eq.13}\E\int_0^{\sigma\wedge\widehat{\tau_y}\wedge\tau}
f(\widehat{S}(t)-\widehat{X}(t))\,
dt=\frac{\kappa(q,0)}{q}\P(\sigma>\tau,\widehat{\tau_y}>\tau)\int_0^\infty
f(u)\Upsilon^q(du)~.\end{equation}
Now we define the \emph{overshoot} at time $\sigma$ as
\begin{equation*}J_\tau:=(\Delta
C(\sigma)-(\widehat{S}(\sigma-)-\widehat{X}(\sigma-)))\cdot
1_{\{\sigma\leq\tau\}}~.\end{equation*}
Using the compensation formula (for details on this formula see for example \cite[p.7]{Ber}) for $f(u):=1_{(z,\infty)}(u)\nu(x+u,\infty)$ and (\ref{2:eq.13}), we get
\begin{align*}&\P(\widehat{S}(\sigma_{-})\leq
y,\widehat{S}(\sigma_{-})-\widehat{X}(\sigma_{-})>z,J_\tau>x,\sigma\leq\tau)\\
&=\E\int_0^{\sigma\wedge\widehat{\tau_y}\wedge\tau}
1_{\{\widehat{S}(t)-\widehat{X}(t)>z\}}\nu(x+\widehat{S}(t)-\widehat{X}(t),\infty)\,dt\\
&=\frac{\kappa(q,0)}{q}\P(\sigma>\tau,\widehat{\tau_y}>\tau)\int_0^\infty
1_{(z,\infty)}(u)\nu(x+u,\infty)\Upsilon^q(du)\\
&=\frac{\kappa(q,0)}{q}\P(\sigma>\tau,\widehat{\tau_y}>\tau)\int_{z+x}^\infty
\nu(u,\infty)\Upsilon^q(du)~.
\end{align*}
When $x\to 0$, $z\to
0$ and $y\to\infty$, we have
\begin{equation*}\P(\sigma\leq\tau)=\frac{\kappa(q,0)}{q}\P(\sigma>\tau)\int_0^{\infty}\nu(u,\infty)\Upsilon^q(du)~.
\end{equation*}
For $z\to 0$ and $y\to\infty$
\begin{equation*}\P(J_\tau>x,\sigma\leq\tau)=\frac{\kappa(q,0)}{q}\P(\sigma>\tau)\int_x^\infty
\nu(u,\infty)\Upsilon^q(du)
\end{equation*} so we have
\begin{equation*}\P(J_\tau>x|\sigma\leq\tau)=\frac{\int_x^\infty
\nu(u,\infty)\Upsilon^q(du)}{\int_0^\infty
\nu(u,\infty)\Upsilon^q(du)}~.
\end{equation*}
These results are summarized in the following lemma.
\begin{lemma}
For $x,~y,~z>0$ we have
\begin{align}&\P(\widehat{S}(\sigma-)\leq
y,\widehat{S}(\sigma-)-\widehat{X}(\sigma-)>z,J_\tau>x,\sigma\leq\tau)\nonumber\\
&=\frac{\kappa(q,0)}{q}\P(\sigma>\tau,\widehat{\tau_y}>\tau)\int_{z+x}^\infty
\nu(u,\infty)\Upsilon^q(du)~,\label{2:eq.14}
\end{align}
\begin{equation}\P(\sigma\leq\tau)=\frac{\kappa(q,0)}{q}\P(\sigma>\tau)\int_0^{\infty}\nu(u,\infty)\Upsilon^q(du)\label{2:eq.15}
\end{equation}
and
\begin{equation}\P(J_\tau>x|\sigma\leq\tau)=\frac{\int_x^\infty
\nu(u,\infty)\Upsilon^q(du)}{\int_0^\infty
\nu(u,\infty)\Upsilon^q(du)}~. \label{2:eq.16}
\end{equation}
\end{lemma}
Now we define the equivalent of the integrated tail of the distribution of the jumps, i.e.
\begin{equation}\label{2:eq.17}H_\tau(x):=\frac{\int_0^x \nu(u,\infty)\Upsilon^q(du)}{\int_0^\infty
\nu(u,\infty)\Upsilon^q(du)}~.
\end{equation}
Then we have
\begin{equation*}
\P(J_\tau>x|\sigma\leq\tau)=1-H_\tau(x)~,~~x>0~.
\end{equation*}
Letting $z\to 0$ and $x\to 0$ in (\ref{2:eq.14}) leads to the fact that $1_{\{\sigma\leq\tau\}}$ and $\widehat{S}((\sigma\wedge\tau)-)$ are independent. More precisely, we have
\begin{align*}\P(\widehat{S}((\sigma\wedge\tau)-)\leq y,\sigma\leq\tau)&=\P(\widehat{S}(\sigma-)\leq y,\sigma\leq\tau)\\
&=\P(\sigma>\tau,\widehat{S}(\tau-)\leq
y)\frac{\kappa(q,0)}{q} I~,\end{align*}
where
\begin{equation*}
I:=\int_0^\infty
\nu(u,\infty)\Upsilon^q(du)~.
\end{equation*}
On the other hand, using (\ref{2:eq.15}), we have
\begin{equation*}\P(\sigma\leq\tau)=\frac{\frac{\kappa(q,0)}{q}
I}{1+\frac{\kappa(q,0)}{q}I}
\end{equation*} and
\begin{align*}\P(\widehat{S}((\sigma\wedge\tau)-)\leq
y)
&=\P(\widehat{S}((\sigma\wedge\tau)-)\leq
y,\sigma\leq\tau)+\P(\widehat{S}((\sigma\wedge\tau)-)\leq
y,\sigma>\tau)\\
&=\P(\sigma>\tau,\widehat{S}(\tau-)\leq
y)\frac{\kappa(q,0)}{q} I+\P(\widehat{S}(\tau-)\leq
y,\sigma>\tau)\\
&=(1+\frac{\kappa(q,0)}{q} I)\P(\widehat{S}(\tau-)\leq
y,\sigma>\tau)~.\end{align*}
So we get
\begin{align*}
\P(\widehat{S}((\sigma\wedge\tau)-)\leq y)\cdot\P(\sigma\leq\tau)
&=(1+ \frac{\kappa(q,0)}{q}I)\P(\widehat{S}(\tau-)\leq
y,\sigma>\tau)\frac{\frac{\kappa(q,0)}{q}
I}{1+\frac{\kappa(q,0)}{q} I}\\
&=\P(\widehat{S}(\tau-)\leq
y,\sigma>\tau)\frac{\kappa(q,0)}{q} I\\
&=\P(\widehat{S}((\sigma\wedge\tau)-)\leq y,\sigma\leq\tau)~.
\end{align*}
This leads to the following result.
\begin{lemma}
The random event $\{\sigma\leq\tau\}$ and the random variable $\widehat{S}((\sigma\wedge\tau){-})$ are independent.
\end{lemma}
Considering our assumptions, we can define times
\begin{equation*}\sigma_1:=\sigma=\inf\{t>0~:~\Delta
C(t)>\widehat{S}(t-)-\widehat{X}(t-)\}\end{equation*} and
\begin{equation*}\sigma_{n+1}=\inf\{t>\sigma_n~:~\Delta
C(t)>\widehat{S}(t-)-\widehat{X}(t-)\}~,~~n\geq 1~,
\end{equation*}
so it is valid that
\begin{equation*}0<\sigma_1<\sigma_2<\ldots~, a.s.
\end{equation*}
Now we can use all of the above results to decompose the process $\widehat{X}$ until the time $\tau$ using the \emph{modified ladder heights},
\begin{align*}&{L_0}^{\tau}:=\widehat{S}((\sigma_1\wedge\tau)-)~,\\
&{J_1}^{\tau}:=\widehat{S}(\sigma_1\wedge\tau)-\widehat{S}((\sigma_1\wedge\tau)-)~,\\
&{L_1}^{\tau}:=\widehat{S}((\sigma_2\wedge\tau)-)-\widehat{S}(\sigma_1\wedge\tau)~~\textrm{on}~~\{\sigma_1\leq\tau\}\end{align*}
and so on, until ${J_{N_\tau}}^{\tau}$ and ${L_{N_\tau}}^{\tau}$, where
$$N_\tau:=\max\{n\in\N~:~\sigma_n\leq\tau\}~.$$
Using the strong Markov property and properties of the exponential distribution, we note that $\P(\sigma_n\leq\tau)={\rho_\tau}^n$, for
$\rho_\tau:=\P(\sigma\leq\tau)$ (the probability which we expressed in the Lemma 2.2.) and that $N_\tau$ has geometric distribution with
parameter $1-\rho_\tau$.\\
\\
Now we can decompose
\begin{equation} \label{2:eq.18}\widehat{S}(\tau)=\sup_{0\leq t\leq \tau}
\widehat{X}(t)={L_0}^{\tau}+{J_1}^{\tau}+{L_1}^{\tau}+\cdots +{J_{N_\tau}}^{\tau}+{L_{N_\tau}}^{\tau}~,
\end{equation}
so
\begin{align*}&\P(\widehat{S}(\tau)\leq
x)=\P({L_0}^{\tau}+{J_1}^{\tau}+{L_1}^{\tau}+\cdots+
{J_{N_\tau}}^{\tau}+{L_{N_\tau}}^{\tau}\leq x)\\
&=\sum_{n=0}^{\infty}
\P({L_0}^{\tau}+{J_1}^{\tau}+{L_1}^{\tau}+\cdots+
{J_{N_\tau}}^{\tau}+{L_{n}}^{\tau}\leq x, N_\tau=n)~.
\end{align*}
Using the independency of the random event $\{\sigma\leq\tau\}$ and random variable $\widehat{S}((\sigma\wedge\tau){-})$ we have
\begin{align*}\P({L_0}^{\tau}\leq x,N_\tau=0)&=\P(\widehat{S}((\sigma\wedge\tau)-)\leq x,\sigma> \tau)\\
&=\P(\widehat{S}((\sigma\wedge\tau)-)\leq x)\cdot\P(\sigma>\tau)=G_\tau(x)\cdot(1-\rho_\tau)~,
\end{align*}
where $G_\tau$ is defined as the distribution function of $\widehat{S}((\sigma\wedge\tau)-)$ and $\rho_\tau=\P(\sigma\leq\tau)$ is as before.
Using the independency from Lemma 2.3., it follows that the conditional distribution of $\widehat{S}((\sigma\wedge\tau)-)$ given $\sigma\leq\tau$
is equal to the unconditional distribution of $\widehat{S}((\sigma\wedge\tau)-)$, so we have
\begin{align*}P({J_1}^{\tau}\leq x,{L_0}^{\tau}\leq
y|\sigma\leq\tau)
&=P({J_1}^{\tau}\leq x|\sigma\leq\tau)\cdot\P({L_0}^{\tau}\leq
y|\sigma\leq\tau)\\
&=H_\tau(x)\cdot G_\tau(y)~.
\end{align*}
Now we have
\begin{align*}\P({J_1}^{\tau}\leq x,{L_0}^{\tau}\leq y,\sigma\leq\tau)
&=P({J_1}^{\tau}\leq x,{L_0}^{\tau}\leq
y|\sigma\leq\tau)\P(\sigma\leq\tau)\\
&=G_\tau(y)(1-H_\tau(x)) \rho_\tau
\end{align*}
and, using the strong Markov property, it follows that
\begin{equation*}\P({L_0}^{\tau}+{J_1}^{\tau}+{L_1}^{\tau}+\ldots
{J_{N_\tau}}^{\tau}+{L_{N_\tau}}^{\tau}\leq
x,N_\tau=n)=(1-\rho_\tau)
{\rho_\tau}^n({G_\tau}^{(n+1)*}*{H_\tau}^{n*})(x)~.
\end{equation*}
This leads us to the main result of the paper. i.e.
\begin{thm}
For the general risk process $X=Y-C$, where $Y$ is a L\'evy process, $C$ an independent subordinator and $\tau$ an independent exponential time ($\tau\sim Exp(q)$, $q>0$), on $[0,\tau]$ under the assumptions (\ref{2:eq.3}) and (\ref{2:eq.4}) we have that
\begin{equation}\label{2:eq.19}\P(\widehat{S}(\tau)\leq
x)=(1-\rho_\tau)\sum_{n=0}^\infty {\rho_\tau}^n(G_\tau^{(n+1)\ast}\ast
H_\tau^{n\ast})(x)~,~~x\geq 0~,\end{equation} for
\begin{equation*}\rho_\tau=\P(\sigma\leq\tau)=\frac{\kappa(q,0)}{q}\P(\sigma>\tau)\int_0^{\infty}\nu(u,\infty)\Upsilon^q(du)~,\end{equation*}
$G_\tau$ the distribution function of $\widehat{S}((\sigma\wedge\tau)-)$ and
\begin{equation*}H_\tau(x)=\frac{\int_0^x \nu(u,\infty)\Upsilon^q(du)}{\int_0^\infty
\nu(u,\infty)\Upsilon^q(du)}~.\end{equation*}
\end{thm}
Now we can apply our result to a more specific case, when $c>0$ is a premium rate, $C$ is a subordinator with L\'evy measure $\nu$ and finite expectation (i.e. $\E[C(1)]<\infty$) which models the claim process (this is a natural choice from the perspective of the risk theory, since the claim process needs to be a nondecreasing process with stationary and independent increments) and the perturbation $Z$ is modelled as a spectrally negative L\'evy process with finite expectation and, without loss of generality, we can assume that $\E[Z(1)]=0$. We also assume that net profit condition is valid, i.e.
\begin{equation}\label{2_eq.20}\E[C(1)]<c~. \end{equation}
Now we observe the generalized risk process
\begin{equation}\label{2:eq.21}
X(t)=ct-C(t)+Z(t)~,~~t\in[0,\tau]~.
\end{equation}
In this setting, it follows that $X$ is the spectrally negative L\'evy process with finite expectation, i.e. $\E[X(1)]<\infty$ and  $\E[X(1)]=c-\E[C(1)]>0$. This means that we can apply the result valid for the supremum of the spectrally negative L\'evy processes, namely that $S(\tau)\sim Exp(\phi(q))$, where $\phi=\phi_X=\psi_X^{-1}$ is, as before, the inverse of the Laplace exponent of the process $X$. Details can be found for example in \cite[Thm VII.1.1]{Ber} and \cite[Cor VII.1.2.]{Ber}.\\
\\
Now in (\ref{2:eq.5}) we have
\begin{align*}
&\E\int_0^{\tau} 1_{\{\widehat{S}(t)-\widehat{X}(t)\leq x\}}\,dt\nonumber\\
&=\frac{1}{q}\int_0^{\infty} qe^{-tq}\P(S(t)\leq x)\,dt\nonumber\\
&=\frac{1}{q}\P(S(\tau)\leq x)=\frac{1}{q}(1-e^{-\phi(q)x})\end{align*}
and we can apply this result to get analogues of the results from Theorem 2.4. following the same procedure as before. On the other hand, directly using the Theorem 2.4. in this specific setting, we get
\begin{align*}
&\P(\sigma\leq\tau)=\P(\sigma>\tau)\frac{\kappa(q,0)}{q}\int_0^\infty \nu(u,\infty)\Upsilon^q(du)\nonumber\\
&=\P(\sigma>\tau)\frac{\kappa(q,0)}{q}\int_0^\infty\Upsilon^q(u)\nu(du)\nonumber\\
&=\P(\sigma>\tau)\frac{1}{q}\int_0^\infty(1-e^{-\phi(q)u})\nu(du)\nonumber\\
&=\P(\sigma>\tau)\int_0^\infty\frac{\phi(q)}{q}e^{-\phi(q)u}\nu(u,\infty)\,du~,
\end{align*}
where the third line follows using (\ref{2:eq.6b}) with the fact that $S(\tau)\sim Exp(\phi(q))$, i.e.
\begin{align*}
&\int_0^\infty e^{-qt}\P(T(x)>t)dt\nonumber\\
&=\frac{1}{q}\int_0^\infty qe^{-qt}\P(S(t)\leq x)\,dt\nonumber\\
&=\frac{1}{q}\P(S(\tau)\leq x)=\frac{1}{q}(1-e^{-\phi(q)x})~.
\end{align*}
Similarly we get
\begin{equation*}\P(J_\tau>x|\sigma\leq\tau)=\frac{1}{\int_0^{\infty}e^{-\phi(q)u}\nu(u,\infty)\,du}
\int_x^{\infty}e^{-\phi(q)(u-x)}\nu(u,\infty)\,du~.
\end{equation*}
Following the same decomposition of the supremum as in the (\ref{2:eq.18}) we have the Pollaczek-Khinchine type formula in this setting.
\begin{corollary} For the risk model described above and $x\geq 0$,
\begin{equation}\label{2:eq.21}\P(\widehat{S}(\tau)\leq x)=(1-p_\tau)\sum_{n=0}^{\infty}
{p_\tau}^n({\widetilde{G}_\tau}^{(n+1)*}*{H_\tau}^{n*})(x)~, \end{equation}
where
\begin{equation*}H_\tau(x)=\frac{1}{\int_0^\infty e^{-\phi(q) u}\nu(u,\infty)
du}\int_0^x e^{-\phi(q) u}\nu(u,\infty)\,du~,
\end{equation*}
$\widetilde{G}_\tau$ is the distribution function of $\widehat{S}((\sigma\wedge\tau)-)$ and \begin{equation*}
p_\tau=\P(\sigma\leq\tau)=\P(\sigma>\tau)\frac{\phi(q)}{q}\int_0^{\infty}
e^{-\phi(q)u}\nu(u,\infty)du~.
\end{equation*}
\end{corollary}

\section{The ladder process}\label{sec-3}
Now we will observe our results and necessary assumptions from the point of view of the ladder process. The explicit formula for the Laplace exponent $\widehat{\kappa}(\alpha,\beta)$ is known only for the spectrally negative case, so we will take a spectrally negative L\'evy process $Y$
and an independent subordinator $C$ and observe the process
\begin{equation}\label{3:eq.1}X(t)=Y(t)-C(t)=ct-C(t)+Z(t)~,~t\geq 0~.\end{equation}
In that case we have an explicit formula for $\kappa$ (see for example \cite[Thm VII.1.4.]{Ber}),
$$\widehat{\kappa}(\alpha,\beta)=k\cdot\frac{\alpha-\psi(\beta)}{\phi(\alpha)-\beta}~,~~\alpha,\beta>0$$
(for constant $k$ we can take $k=1$ without loss of generality).\\
\\
Let us take the following assumptions:
\begin{itemize}\item[(i)] $C$ and $Z$ have finite expectations (and without loss of generality we can take $\E Z(1)=0$),
\item[(ii)] $\E C(1)>c$,
\end{itemize}
 i.e. we take the opposite to the standard net profit condition. Let us notice that we now consider the same setting as in the Corollary 2.5., but now we are not assuming the net profit condition.
We can now see that these assumptions imply that $\E X(1)=c-\E C(1)+\E Z(1)=c-\E C(1)<0$ so
$\widehat{X}\to +\infty$. In other words, if we observe this from the point of view of the ladder process $\widehat{H}$ associated to $\widehat{X}$, this process will be a subordinator, but killed at rate zero, i.e. a non-killed subordinator.\\
\\
Now we have that $\psi_X'(0_+)=\E X(1)<0$ and the Laplace exponent of $X$, $\psi_X$, has no unique root. This implies that $\psi_X(0)=0$ and $\psi_X(b)=0$ for some $b>0$, so the Laplace exponent of $X$ is actually a function $\psi_X:[b,+\infty\rangle\to[0,\infty\rangle$ and we will denote its inverse as $\phi_X$. The Laplace exponent of the ladder process can be written as
\begin{align*}\widehat{\kappa}(\beta)&=\lim_{\alpha\to 0}
\widehat{\kappa}(\alpha,\beta)=\frac{\psi_X(\beta)}{\beta-\phi_X(0)}
=\frac{c\beta-\psi_C(\beta)+\psi_Z(\beta)}{\beta - \phi_X(0)}\\
&=(1+\frac{\phi_X(0)}{\beta - \phi_X(0)})(c-\frac{\psi_C(\beta)}{\beta}+
\frac{\psi_Z(\beta)}{\beta})~.
\end{align*}
For example, if $Z$ is modelled by the Brownian motion (so its Laplace exponent is $\psi_Z(\beta)=\beta^2$, $\beta>0$) and $C$ is the Poisson process (its Laplace exponent is then given by
$\psi_C(\beta)=\lambda\cdot(1-e^{-\beta})$, $\lambda>0$), the Laplace exponent of $X$ is given by
\begin{equation*}\psi_X(\beta)=c\beta-\lambda\cdot(1-e^{-\beta})+\beta^2~.\end{equation*}
If we subtract the drift, i.e. $\frac{\psi_Z(\beta)}{\beta}=
\frac{\beta^2}{\beta}=\beta$, from $\widehat{\kappa}$, we get
\begin{align*}\widehat{\kappa}(\beta)-\beta&=\frac{c\beta-\lambda(1-e^{-\beta})+\beta^2}{\beta-\phi(0)}-\beta\\
&=\frac{(c+\phi_X(0))\beta-\lambda(1-e^{-\beta})}{\beta-\phi(0)}=:\varphi(\beta)~.
\end{align*}
For $\beta\to\infty$, we have
\begin{equation*}
\lim_{\beta\to\infty}
\varphi(\beta)=\lim_{\beta\to\infty}\frac{(c+\phi_X(0))\beta-\lambda(1-e^{-\beta})}{\beta-\phi_X(0)}=
\lim_{\beta\to\infty}
\frac{c+\phi_X(0)-\lambda\frac{1-e^{-\beta}}{\beta}}{1-\frac{\phi_X(0)}{\beta}}=c+\phi_X(0)~,
\end{equation*}
which is obviously finite. But finiteness of the Laplace exponent at infinity means that the L\'{e}vy measure is finite as well and that is valid only in the case of the compound Poisson process. This fact implies that our modified Laplace exponent $\varphi(\beta)$ is equal to the Laplace exponent of the compound Poisson process, i.e. $\varphi(\beta)=\widetilde{\lambda}\int_{(0,\infty)}
(1-e^{-\beta x})\,\widetilde{F}(dx)$,
where $\widetilde{\lambda}$ is the intensity of the jumps and
$\widetilde{F}$ distribution function of the jumps of the compound Poisson process. Then we have
\begin{equation*}\lim_{\beta\to\infty} \varphi(\beta)=\lim_{\beta\to\infty}
\widetilde{\lambda}\int_{(0,\infty)} (1-e^{-\beta
x})\widetilde{F}(dx)=\widetilde{\lambda}~,\end{equation*}
which means that the jumps of the subordinator $\widehat{H}$  behave just like the compound Poisson process with the intensity of the jumps equal to $c+\phi_X(0)$.\\
\\
For some general subordinator $C$ we have
$$\psi_C(\beta)=\int_{(0,\infty)} (1-e^{-\beta x})\,\nu(dx)~,$$
so again
\begin{equation*}\varphi(\beta)=\widehat{\kappa}(\beta)-\beta=\frac{(c+\phi_X(0))\beta-\int_{(0,\infty)}
(1-e^{-\beta x})\nu(dx)}{\beta -\phi_X(0)}
\end{equation*}
and
\begin{equation*}\lim_{\beta\to\infty} \varphi(\beta)=\lim_{\beta\to\infty}
\frac{c+\phi_X(0)-\int_{(0,\infty)}\frac{1-e^{-\beta
x}}{\beta}}{1-\frac{\phi_X(0)}{\beta}}=c+\phi_X(0)~,
\end{equation*}
where we used the fact that the functions under the integral are bounded with $1-e^{-x}\leq x$ and
$\int_{(0,\infty)} x\nu(dx)=\E C(1)<\infty$, as we assumed. If we would take some other perturbation $Z$ instead of the Brownian motion, we would again obtain the same result in $\widehat{\kappa}$.\\
\\
Let us notice that for some generalized risk process $X'$ in the standard setting, i.e.
$$X'=Y-D~,$$
where $D$ is a subordinator with finite expectation and $\E D(1)<c$ (so, for $X'$ the net profit condition is valid) and $\widehat{\kappa'}$ the Laplace exponent of the ladder height process associated to the process $X'$ we would have:
\begin{equation*}\widehat{\kappa'}(\beta)=c-\frac{\psi_C(\beta)}{\beta}+\frac{\psi_Z(\beta)}{\beta}~.
\end{equation*}
So in this case we have
\begin{equation*}\widehat{\kappa'}(\beta)-\frac{\psi_Z(\beta)}{\beta}=c-\frac{\psi_C(\beta)}{\beta}=c-\int_{0}^{\infty}
e^{-\beta x}\nu(x,\infty)\,dx~.
\end{equation*}
In other words, in this case we have that
$\widehat{\kappa}(0+)\geq 0$ if and only if
$c-\E C(1)=c-\int_0^{\infty} \nu(x,\infty)\,dx\geq 0$
if and only if $\E C(1)\leq c$ and $\lim_{\beta\to\infty} (c-\int_{0}^{\infty}
e^{-\beta x}\nu(x,\infty)\,dx)$ is always equal to $c$.\\
\\
What does this mean? It means that from the point of view of the ladder process, the net profit condition doesn't play the important role in decomposing the supremum of the observed process - the only factor that changes is the value $\phi_X(0)$ which is not equal to zero in the case that net profit condition is not satisfied and this value affects the rate of the underlying ladder process. This observation is consistent with the results we obtained in the Theorem 2.4. and Corollary 2.5. More precisely, dropping the standard net profit condition still allows us to decompose supremum of the dual of the generalized risk process at modified ladder epochs and achive a Pollaczek-Khinchine type results and the key quantity that makes the change in the results is the inverse of the Laplace exponent of the generalized risk process which arises from the ladder height process in the background. If we also drop the assumptions on the finite expectation, the same type of decomposition is obtainable (under some necessary assumptions as we have seen in the Section 2), only in terms of the more general renewal functions, as we have shown in the Theorem 2.4.\\
\\
\textbf{Acknowledgement:} This work has been supported in part by Croatian Science Foundation under
the project 3526.

\end{document}